%Paper: alg-geom/9505010
%From: lbadescu@ROIMAR.IMAR.RO
%Date: Mon, 08 May 1995 11:28:02 0200

\input amstex
\documentstyle{amsppt}
\NoBlackBoxes
\pagewidth{29pc}
\pageheight{40pc}
\topmatter
\def\Lef{\operatorname{Lef}}

\def\Proj{\operatorname{Proj}}
\def\Spec{\operatorname{Spec}}
\def\codim{\operatorname{codim}}

\def\Ass{\operatorname{Ass}}
\def\Spf{\operatorname{Spf}}
\def\deg{\operatorname{deg}}

\title{ALGEBRAIC BARTH-LEFSCHETZ THEOREMS}\endtitle
\rightheadtext{ALGEBRAIC BARTH-LEFSCHETZ THEOREMS}
\leftheadtext{LUCIAN B\v{A}DESCU}
\author{LUCIAN B\v{A}DESCU}\endauthor
\address{Institute of Mathematics of the Romanian Academy,
P.O. Box 1-764, Ro-70700 Bucharest, Romania}\endaddress
\email{lbadescu\@ imar.ro}\endemail
%\keywords{ }\endkeywords
%\abstract{ }\endabstract
%\thanks{This paper is in final form and no version of
%it will be submitted elsewhere}\endthanks
\endtopmatter
\document

\head{0. Introduction}\endhead

We shall work over a fixed algebraically closed field $k$ of
arbitrary characteristic. By an algebraic variety over $k$ we shall
mean a reduced algebraic scheme over $k$. Fix a positive integer $n$ and
$e=(e_0,e_1,...,e_n)$ a system of $n+1$ weights (i.e. $n+1$ positive
integers $e_0,e_1,...,e_n$). If $k[T_0,T_1,...,T_n]$ is the polynomial
$k$-algebra in $n+1$ variables, graded by the conditions
${\deg}(T_i)=e_i, i=0,1,...,n$, denote by ${\Bbb P}^{n}(e)=
{\Proj}(k[T_0,T_1,...,T_n])$ the $n$-dimensional weighted
projective space over $k$ of weights $e$. We refer the reader to [3]
for the basic properties of weighted projective spaces.
According to Zariski ([22], see also [16], [15]), if $Y$ is a closed subscheme
of an algebraic variety $X$, one can define the ring $K(\hat{X}_{/Y})$ of
formal rational functions of $X$ along $Y$. Then $K(\hat{X}_{/Y})$
is a $k$-algebra, and there is a canonical map of $k$-algebras
$K(X)\rightarrow K(\hat{X}_{/Y})$, where $K(X)$ is the usual ring
of rational functions of $X$ ($K(X)$ is a field if $X$ is irreducible).
According to [16], $Y$ is said to be $\bold{G}_3$ in $X$
if this map is an isomorphism. Let $X$ be an arbitrary algebraic
scheme over $k$, and let $d\geq 0$ be a non-negative integer. Then
$X$ is said to be $d$-connected if every irreducible component of
$X$ is of dimension $\geq d+1$ and if $X\setminus W$ is connected
for every closed subscheme $W$ of $X$ of dimension $<d$. For example,
$X$ is $0$-connected if $X$ is connected and of dimension $\geq 1$; an
irreducible algebraic variety $X$ of dimension $n\geq 1$ is always
$(n-1)$-connected.

Then the main result of this paper is the following.

\proclaim{Theorem (0.1)}
Let $f:X\rightarrow {\Bbb P}^{n}(e)\times{\Bbb P}^{n}(e)$
be a finite morphism from a $d$-connected algebraic variety $X$ such that
$d\geq n$. Then $f^{-1}(\Delta)$ is $(d-n)$-connected, where $\Delta$ is
the diagonal of ${\Bbb P}^{n}(e)\times{\Bbb P}^{n}(e)$. Moreover,
$f^{-1}(\Delta)\setminus W$ is $\bold{G}_3$ in $X\setminus W$ for every
closed subscheme $W$ of $f^{-1}(\Delta)$ of dimension $<d-n$.
\endproclaim

In the case of ordinary projective spaces (i.e. when $e_i=1$ for every
$i=0,1,...,n$) the first statement of theorem 1 is well known in the
literature as the Fulton-Hansen connectedness theorem (see [5], or also
[6], or [12]). As far as I know the last part of theorem 1 is new even
for ordinary projective spaces. This latter fact has some
interesting consequences, namely the following two theorems.

\pageheight{47pc}
\proclaim{Theorem (0.2)} Let $Y$ be a closed irreducible subvariety of
${\Bbb P}^{n}(e)$ of dimension $>\frac{n}{2}$. Then ${\Delta}_Y\setminus W$
is $\bold{G}_3$ in $Y\times Y\setminus W$ for every closed subscheme
$W$ of the diagonal ${\Delta}_Y$ of $Y\times Y$ such that
${\dim}(W)<2{\dim}(Y)-n-1$. If moreover $Y$ has (locally) the ${\bold S}_2$
property of Serre (e.g. if $Y$ is normal), then for every vector bundle
$E$ on $Y\times Y$ the natural map
$$ H^{0}(Y\times Y,E)\rightarrow H^{0}(\widehat{Y\times Y},\hat{E}) $$
is an isomorphism, where $\hat{E}$ is the formal completion of $E$
along ${\Delta}_Y$, and
$\widehat{Y\times Y}=\widehat{Y\times Y}_{/{\Delta}_Y}$. In other words,
the weak Grothendieck-Lefschetz condition ${\Lef}(Y\times Y,{\Delta}_Y)$
holds (see [9]).
\endproclaim

An immediate consequence of theorem (0.2) via a result of Speiser (see
[21], or also theorem (1.11) below) is the following.

\proclaim{Theorem (0.3)} Let $Y$ be a closed irreducible subvariety of
${\Bbb P}^{n}(e)$ of dimension $>\frac{n}{2}$ which is (locally)
${\bold S}_2$.Then every stratified vector bundle on $Y$ is trivial.
\endproclaim

For the definition of stratified vector bundles see [10], or also [20].
This definition and the basic properties of stratified vector bundles
will also be briefly recalled in the first section.
In some special cases theorem 3 and the last part of theorem 2 were
known before. Specifically, if $Y$ is a local complete intersection
in ${\Bbb P}^n$ over a field $k$ of characteristic zero, theorem 3
and the last part of theorem 2 were proved by Ogus in [18]. Note
that when $\text{char}(k)=0$ the concept of stratified vector bundle
is the same as the one of vector bundle with integrable connection
(see [18]). If $\text{char}(k)>0$ and $Y$ is an irreducible locally
Cohen-Macaulay of ${\Bbb P}^n$ these results were proved by Speiser
in [21]. Note that the methods of Ogus (in characteristic $0$) or those of
Speiser (in positive characteristic) are completely different from the
methods used in this paper. Our approach (which is based on results of
Hironaka-Matsumura [16] and of Faltings [4] involving formal rational
functions) offer therefore not only characteristic free proofs but
also more general results.

A first version of this paper was written during my visit at the
Universities of Pisa (December 1993) and Ferrara (January 1994). I
am grateful to Fabrizio Catanese and Alexandru Lascu for their kind
invitation and excellent atmosphere I found there. I also want to
thank G. Chiriacescu who read some parts of the paper and made
useful remarks and suggestions.

\head{1. Background material}\endhead

In this section we gather together some known results we are going
to use in this paper. Let $X$ be an algebraic variety, and let
$X_1,...,X_m$ be the irreducible components of $X$. Fix $d\geq 0$ a
non-negative integer.

\definition{Definition (1.1)} ([13, [12]) A sequence $Z_0,Z_1,...,Z_n$ of
(not necessarily mutually distinct) irreducible components of $X$ is called
a $d$-join within $X$ if ${\dim}(Z_i)\geq d+1$ for every $i=0,1,...,n$
and if ${\dim}(Z_{j-1}\cap Z_j)\geq d$ for every $j=2,...,n$.
\enddefinition

The following elementary fact will be useful.

\proclaim{Proposition (1.2)} (Hartshorne [13], or also [12])
An algebraic variety $X$ is $d$-connected if and only if
$X=Z_0\cup Z_1\cup ...\cup Z_n$ for some $d$-join
$Z_0,Z_1,...,Z_n$ within $X$.
\endproclaim

\proclaim {Theorem (1.3)} (Grothendieck) Let $X$ be a
$d$-connected algebraic variety over $k$, and let
$f:X\rightarrow {\Proj}(S)$ be a finite morphism, where $S$
is a finitely generated graded $k$-algebra. Let $t_1,...,t_r
\in S_{+}$ be homogeneous elements of positive degrees. If $d\geq r$ then
$f^{-1}(V_+(t_1,...,t_r))$ is $(d-r)$-connected. Moreover, if $X$ is
irreducible and ${\dim}(X)\geq r$ then $f^{-1}(V_{+}(t_1,...,t_r))$ is
non-empty.
\endproclaim

\remark {Remark} Theorem (1.3) can be found (in a slightly different
formulation) in [9],$\acute{e}$xpos$\acute{e}$ XIII. In the appendix
another proof based on the so-called Hartshorne-Lichtenbaum theorem
can be found.
\endremark

\proclaim{Theorem (1.4)} (Hironaka-Matsumura [16])
Let $f:X'\rightarrow X$ be a proper surjective morphism of algebraic
varieties, with $X$ irreducible and such that every irreducible
component of $X'$ dominates $X$. Let $Y$ be a closed subscheme of
$X$ and set $Y':=f^{-1}(Y)$. Then there is a canonical isomorphism
of $k$-algebras
$$ K(\hat{X'}_{/Y'})\cong [K(\hat{X}_{/Y})\otimes _{K(X)}K(X')]_0 ,$$
where $K(X')$ is the ring of rational functions of $X'$,
$K(\hat{X'}_{/Y'})$ is the ring of formal rational functions of $X'$
along $Y'$, and $[A]_0$ denotes the total ring of fractions of a
commutative ring $A$.
\endproclaim

The next result will play a crucial role in the proof of theorem (0.1).

\proclaim{Theorem (1.5)} (Faltings [4]) Let $X$ be an irreducible closed
subvariety of ${\Bbb P}^n$, and let $L$ be a linear subspace of
${\Bbb P}^n$ of codimension $r$ such that $r<{\dim}(X)$. Then
$X\cap L\setminus W$ is ${\bold G}_3$ in $X\setminus W$ for every closed
subvariety $W$ of $X\cap L$ such that ${\dim}(W)<{\dim}(X)-r-1$
(${\dim}(W)=-1$ if $W=\emptyset$).
\endproclaim

\remark{Remark} Theorem (1.5) is a special case of a result of Faltings
[4]. In the case when $X$ and $L$ have a proper intersection, this
result was proved earlier by Hironaka and Matsumura [16] (in which
case the proof is much more elementary).
\endremark

\proclaim{Proposition (1.6)} ([16]) Let $X$ be an irreducible algebraic
variety, and let $Y$ be a closed subscheme of $X$. Let $u:X'\rightarrow X$
be the normalization morphism. Then $K(\hat{X}_{/Y})$ is a field if and
only if $u^{-1}(Y)$ is connected.
\endproclaim

The next proposition is well known for the rings of usual rational
functions, and should be known in general, but we have no reference
for it (therefore we include a proof).

\proclaim{Proposition (1.7)} Let $X$ be a quasi-projective variety having the
irreducible components $X_1,...,X_m$ (with the reduced structure), and
let $Y$ be a closed subscheme of $X$, such that $Y_i:=Y\cap X_i\neq\emptyset$
for every $i=1,...,m$. Then there is a canonical isomorphism of $k$-algebras
$$ K(\hat{X}_{/Y})\cong K(\hat{X_1}_{/Y_1})\times ...\times
K(\hat{X_m}_{/Y_m}). $$
\endproclaim

\demo{Proof} We shall first prove proposition (1.7) in the case when
$X={\Spec}(A)$ is affine, with $A$ a reduced finitely generated $k$-algebra.
Set ${\Ass}(A)=\lbrace p_1,...,p_m \rbrace$, where every $p_i$ is a
minimal prime ideal of $A$ such that $X_i=V(p_i)$, $i=1,...,m$ and
$$ p_1\cap ...\cap p_m = (0). \tag 1.7.1 $$
Then $Y=V(I)$, with $I$ an ideal of $A$. The hypotheses imply that
$I_i:=I(A/p_i)\neq A/p_i$ for every $i=1,...,m$. Let $\hat{A}$ be the $I$-adic
completion of $A$. Then $\hat{X}_{/Y}={\Spf}(\hat{A})$, and by [16],
(1.1) we have
$$ K(\hat{X}_{/Y})=[\hat{A}]_0, \tag 1.7.2 $$
where, as above, $[\hat{A}]_0$ denotes the total ring of fractions of
$\hat{A}$. Set $q_i=p_i\hat{A}$, $i=1,...m$. Taking into account that
the homomorphism $A\rightarrow \hat{A}$ if flat, (1.7.1) yields
$$ q_1\cap ...\cap q_m=0 ,$$
which implies
$$ q_1[\hat{A}]_0\cap ...\cap q_m[\hat{A}]_0 = (0). \tag 1.7.3 $$
If we denote by $S$ (resp. by $\hat{S}$) the multiplicative system of
all non-zero divizors of $A$ (resp. of $\hat{A}$), then
$p_iS^{-1}A+p_jS^{-1}A=S^{-1}A$ for every $i\neq j$, which implies
$p_iS^{-1}\hat{A}+p_jS^{-1}\hat{A}=S^{-1}\hat{A}$ for every $i\neq j$,
or else, $q_iS^{-1}\hat{A}+q_jS^{-1}\hat{A}=S^{-1}\hat{A}$ for every
$i\neq j$. Since $S\subset \hat{S}$ (because $A\rightarrow \hat{A}$
is flat), we therefore get:
$$ q_i[\hat{A}]_0 + q_j[\hat{A}]_0 = [\hat{A}]_0 \;\text{for}\; i\neq j.
\tag 1.7.4 $$
Then (1.7.2), (1.7.3) and (1.7.4) together with the Chinese Remainder
Theorem imply
$$ K(\hat{X}_{/Y})\cong [\hat{A}]_0/q_1[\hat{A}]_0\times ...\times
[\hat{A}]_0/q_m[\hat{A}]_0.  \tag 1.7.5 $$
To prove the proposition in case $X$ is affine it will be sufficient
to show that
$$ [\hat{A}]_0/q_i[\hat{A}]_0 \cong K(\hat{X_i}_{/Y_i}) \;\text{for
every}\; i=1,...,m.  \tag 1.7.6 $$
Since $K(\hat{X_i}_{/Y_i})\cong [\hat{A}/q_i]_0$ ([16], (1.1)),
then (1.7.6) is equivalent to
$$ {\hat S}^{-1}(\hat{A}/q_i)=[\hat{A}/q_i]_0 \; \text{for every} \;
i=1,...,m.  \tag 1.7.7 $$

To prove (1.7.7), we first show that for every fixed $i$ the ring
$\hat{A}/q_i$ is reduced. To see this, observe that $\hat{A}/q_i$
is the $I$-adic completion of $A/p_i$. So, we have to prove that if
$B$ is a reduced, finitely generated $k$-algebra, and $J\neq B$ is an ideal
of $B$ then the $J$-adic completion $B^{\ast}$ of $B$ is also reduced.
To show this, it will be sufficient to check that for every maximal
ideal $m^{\ast}$ of $B^{\ast}$, the localization ${B^{\ast}}_{m^{\ast}}$
is reduced. Since $B^{\ast}$ is $JB^{\ast}$-adically complete,
$JB^{\ast}\subset m^{\ast}$. Since $B^{\ast}/J^{n}B^{\ast}\cong B/J^n$,
it follows that for every $n\geq 1,
{B^{\ast}}_{m^{\ast}}/J^{n}{B^{\ast}}_{m^{\ast}}\cong B_m/J^{n}B_m$,
where $m:=m^{\ast}\cap B$. From the latter isomorphism we infer
that for every $n\geq 1,
{B^{\ast}}_{m^{\ast}}/m^{\ast n}{B^{\ast}}_{m^{\ast}}\cong B_m/m^{n}B_m$.
Passing to the inverse limits, we get
$\widehat{B^{\ast}}_{m^{\ast}}\cong \hat B_m$
(the completions with respect to the maximal ideals $m^{\ast}$
and $m$ respectively). Now, use a theorem of Chevalley (see [19], IV-4) to
show that the latter ring is reduced. Finally, since
$\widehat{B^{\ast}}_{m^{\ast}}$ is reduced, ${B^{\ast}}_{m^{\ast}}$ is
also reduced, and we are through.

So, we have seen that $\hat{A}/q_i$ is reduced. In particular, we get
$$ q_i=p_{i1}\cap ...\cap p_{in_i}, $$
where $\text{Ass}(\hat{A}/q_i)=\lbrace p_{i1},...,p_{in_i}\rbrace$.

Then by [19], IV-4, ${\Ass}(\hat{A})=\lbrace p_{11},...,p_{1n_1},...,
p_{m1},...,p_{mn_m}\rbrace $. Moreover, all the
$p_{ij}$ 's are minimal prime ideals of $\hat{A}$.

The inclusion ${\hat S}^{-1}(\hat{A}/q_i)\subset
[\hat{A}/q_i]_0$ being clear, it will be sufficient to show that
every non-zero divisor $\hat{a}:=a$ (mod $q_i$) of $\hat{A}/q_i$ is the
class (modulo $q_i$) of an element of $\hat{S}$ (i.e. $\hat{a}=\hat{b}$
with $b$ a non-zero divisor of $\hat{A}$, or else, $b\notin p_{ij}$
for all $i$ and $j$). Take for example $i=1$.
Since $\hat{a}$ is a non-zero divisor in $\hat{A}/q_1$,
$a\notin\lbrace p_{11}\cup ...\cup p_{1n_1}\rbrace$.
Set $\Lambda :=\lbrace (i,j)/a\in p_{ij}\rbrace$. If
$\Lambda$ is empty we have nothing to prove. Suppose therefore
$\Lambda\neq\emptyset$. Then set $J:={\bigcap}_{(i,j)\notin\Lambda} p_{ij}$
and $J':={\bigcup}_{(i,j)\in\Lambda}p_{ij}$. Since all the ideals
$p_{ij}$ are prime, by elementary general facts about prime ideals
(see e.g. [19], I-2), there is an element $u\in J\setminus J'$.
In other words, $u\in p_{ij}$ if and only if $a\notin p_{ij}$.
This shows that $a+u\notin p_{ij}$ for every $i=1,...,m, j=1,...,n_i$.
On the other hand, since $a\notin p_{1j}$ for every $j=1,...,n_1$,
whence $u\in p_{11}\cap ...\cap p_{1n_1}=q_1$. Therefore $\hat{a}=
\hat{b}$, with $b=a+u$, and $b$ not a zero divisor in $\hat{A}$.
Therefore proposition (1.7) is proved if $X$ is affine.

If $X$ is not affine, fix a finite set $A$ of closed points of $Y$
such that $A\cap Y_i\neq\emptyset$ for every $i=1,...,m$. Since $X$ is
quasi-projective, we can find an affine cover
${\lbrace U_{\alpha}\rbrace}_{\alpha}$ of $X$ such that for every
$\alpha$, $A\subset U_{\alpha}$. Since for every $\alpha$ and $\beta$,
$U_{\alpha\beta}:=U_{\alpha}\cap U_{\beta}$ is also
affine and contains $A$, we know the proposition for $U_{\alpha}$ and
for $U_{\alpha\beta}$.
Then everything follows from the statement in the affine case
(already proved) and from the exact diagram (see [16] for details):
$$ \CD
K(\hat{X}_{/Y})\rightarrow
{\prod}_{\alpha} K(\hat{U_{\alpha}}_{/Y\cap U_{\alpha}})
{\rightarrow\atop\rightarrow}
{\prod}_{\alpha ,\beta} K(\hat{U_{\alpha\beta}}_{/Y\cap U_{\alpha\beta}}),
\endCD $$
which reduces the verification to the affine case. Q.E.D.
\enddemo

\proclaim{Corollary (1.8)} In the hypotheses of proposition (1.7),
if $Y_i$ is ${\bold G}_3$ in $X_i$ for every $i=1,...,m$, then
$Y$ is also ${\bold G}_3$ in $X$.
\endproclaim

\demo{Proof}. Direct consequence of proposition (1.7) and of
the well known isomorphism
$$ K(X)=K(X_1)\times ...\times K(X_m).$$
Q.E.D.
\enddemo

Corollary (1.8) allows one to generalize the result of Faltings
stated above (theorem (1.5)):

\proclaim{Corollary (1.9)} Let $X$ be a $d$-connected closed
subvariety of ${\Bbb P}^n$, and let $L$ be a linear subspace of
${\Bbb P}^n$ of codimension $r$ such that $r\leq d$. Then
$X\cap L\setminus W$ is ${\bold G}_3$ in $X\setminus W$ for
every closed subscheme $W$ of $X\cap L$ such that
${\dim}(W)<d-r$.
\endproclaim

\demo{Proof} Let $X_1,...,X_m$ be the irreducible components of $X$.
Since $X$ is $d$-connected, ${\dim}(X_i)\geq d+1$, and in particular,
$r<{\dim}(X_i)$ for every $i=1,...,m$. By theorem (1.5),
$X_i\cap L\setminus W$ is ${\bold G}_3$ in $X_i\setminus W$, i.e.
$K(X_i\setminus W))\cong K(\widehat{X_i\setminus W}_{/X_i\cap L\setminus W})$
for every $i=1,...,m$ and for every closed subscheme $W$ of $X\cap L$
such that ${\dim}(W)<d-r$. On the other hand, since the irreducible
components of $X\setminus W$ are $X_i\setminus W, i=1,...,m$  everything
follows from corollary (1.8). Q.E.D.
\enddemo

(1.10) Now we recall briefly some definitions and basic facts from
Grothendieck's theory of stratified vector bundles and the descent
theory of faithfully flat morphisms (see [10], [8], or also [21]).
Let $Y$ be an algebraic variety over $k$. Consider the products and
projections:
$$\CD   Y @<p_1<< Y\times Y @>p_2>> Y  \endCD $$
and
$$\CD
Y\times Y @<p_{21}<< Y\times Y\times Y @>p_{31}>> Y\times Y \\
                 @.  @VVp_{32}V  \\
                 @.  Y\times Y
\endCD $$

Denote by $\Delta={\Delta}_Y$ the diagonal of $Y\times Y$, and let
${\Delta}'$ be the diagonal subscheme of $Y\times Y\times Y$. For
every $r\geq 0$ denote by ${\Delta}_r$ (resp. ${\Delta}'_r$) the
$r$-th infinitesimal neighbourhood of $\Delta$ in $Y\times Y$
(resp. of ${\Delta}'$). Then from the previous diagrams we get
projections
$$\CD Y @<q^r_1<< {\Delta}_r @>q^r_2>> Y  \endCD $$
and
$$\CD
{\Delta}_r @<q^r_{21}<< {\Delta}'_r @>q^r_{31}>> {\Delta}_r \\
                     @. @VVq^r_{32}V  \\
		     @. {\Delta}_r
\endCD $$

Let $F$ be be a vector bundle on $Y$. By descent data on $F$ we mean
an isomorphism
$$ \varphi:p^{\ast}_1 (F)\rightarrow p^{\ast}_2 (F) $$
over $Y\times Y$ such that the cocycle condition
$$
p^{\ast}_{31}(\varphi)=p^{\ast}_{32}(\varphi)\circ p^{\ast}_{21}(\varphi)
$$
holds over $Y\times Y\times Y$.

{}From Grothendieck's faithfully flat descent theory (see [8], VIII) applied
to the structural morphism $p:Y\rightarrow {\Spec}(k)$ we know that
descent data hold for $F$ if and only if $F$ is the pull-back via
$p$ of a vector bundle over ${\Spec}(k)$, i.e. if and only if
$F$ is a trivial vector bundle. This result is a very useful general
criterion for a vector bundle to be trivial.

A stratification on $F$ is a compatible system of isomorphisms
$$ {\varphi}_r:(q^r_1)^{\ast}(F)\rightarrow (q^r_2)^{\ast}(F) $$
over ${\Delta}_r$ for all $r\geq 0$ such that ${\varphi}_0=\text{id}$
and the cocycle condition
$$ (q^r_{31})^{\ast}({\varphi}_r)=
(q^r_{32})^{\ast}({\varphi}_r)\circ (q^r_{21})^{\ast}({\varphi}_r) $$
holds on ${\Delta}'_r$.

A stratified vector bundle on $Y$ is a vector bundle $F$ with a
stratification on it. In other words, giving a stratification on $F$
is the same as giving "formal descent data" on $F$. That is to say,
if $X$ is the formal completion of $Y\times Y$ along $\Delta$,
and if $X'$ is the formal completion of $Y\times Y\times Y$
along ${\Delta}'$, then we get projections
$$\CD   Y @<q_1<< X @>q_2>> Y  \endCD $$
and
$$\CD
        X @<q_{21}<< X' @>q_{31}>> X  \\
	         @.  @VVq_{32}V \\
                 @.  X
\endCD $$
and "formal descent data" on $F$ consist of an isomorphism
$$  \psi:q^{\ast}_1 (F)\rightarrow q^{\ast}_2 (F)  $$
on $X$ together with the cocycle condition
$$ q^{\ast}_{31}(\psi)=q^{\ast}_{32}(\psi)\circ q^{\ast}_{21}(\psi) $$
on $X'$. In characteristic zero giving a stratification on
$F$ is the same as giving an integrable connection on $F$ (see [10],
[18]). The approach of vector bundles endowed with integrable
connections and the De Rham cohomology associated to them was
taken up by Ogus to give algebraic proofs of Barth's theorems
(see [18]). On the other hand,Gieseker and Speiser provided further
methods to study stratifications in positive characteristic
(see [7], [21]).The relationship between stratified vector bundles
and descent data is given by the following result of Speiser:
\proclaim{Theorem (1.11)} (Speiser) Let $Y$ be an irreducible
algebraic variety over $k$ for which the diagonal $\Delta$ of
$Y\times Y$ satisfies the weak Grothendieck-Lefschetz condition
${\Lef}(Y\times Y,\Delta)$. Then giving a stratification on a vector
bundle $F$ on $Y$ is equivalent to giving descent data on $F$.
In particular, if ${\Lef}(Y\times Y,\Delta)$ holds, then every
stratified vector bundle on $Y$ is trivial.
\endproclaim

In fact Speiser proved in [21] the first part of theorem (1.9);
the last part follows from the first one and from
Grothendieck's faithfully flat descent theory [8], VIII.

\head{2. Proof of theorem (0.1)} \endhead

To prove the connectivity part of theorem (0.1) we use the same
main ideas of the proof of Fulton-Hansen connectivity theorem given
in [3] or in [6]. First we shall need the following more general
version of theorem (1.3) of Grothendieck.

\proclaim{Proposition (2.1)} Let $S$ be a finitely generated graded
$k$-algebra, $t_1,...,t_r\in S_{+}$ homogeneous elements of
positive degrees, and $U$ a Zariski open subset of ${\Proj}(S)$
containing $L:=V_{+}(t_1,...,t_r)$. Let $f:X\rightarrow U$ be a
finite morphism, with $X$ a $d$-connected algebraic variety over
$k$. If $d\geq r$ then $f^{-1}(L)$ is $(d-r)$-connected.
\endproclaim

\demo{Proof} First we shall prove the proposition in the case when
$X$ is irreducible. Then $X$ is $({\dim}(X)-1)$-connected. Therefore
in this case the hypothesis reads ${\dim}(X)>r$. By passing to the
normalization we may assume that $X$ is normal. Let $Z'$ be the
closure of $X':=f(X)$ in $P:={\Proj}(S)$, and let $g:Z\rightarrow Z'$
be the normalization of $Z'$ in the field $K(X)$ of rational functions
of $X$ (which makes sense because the dominant morphism $X\rightarrow Z'$
yields the finite field extension $K(Z')=K(X')\rightarrow K(Z)$). Then
we get a commutative diagram of the form
$$\CD        X  @>i>>  Z \\
             @VfVV     @VVgV \\
	     X' @>>i'> Z'
\endCD $$
in which $i$ and $i'$ are open immersions ($i$ is an open immersion
because $X$ is normal), and $g$ is a finite morphism. Since $L\subset U$
and $X"\subset U$, then $X'\cap L=Z'\cap L$, whence $f^{-1}(L)=g^{-1}(L)$.
Then proposition (2.1) follows (in case when $X$ is irreducible) from
theorem (1.3) applied to the composition of the finite morphism
$g:Z\rightarrow Z'$ followed by the closed immersion $Z'\subset P$.

Assume now $X$ reducible. In this case we shall first prove the proposition
in the special case $d=r$. Then we simply have to show that
$Y:=f^{-1}(L)$ is connected and non-empty. Let $X_1,...,X_m$ be the
irreducible components of $X$ (all of dimension $\geq d+1=r+1$ since $X$
is $d$-connected). If we set $f_i:=f/X_i$ then by what we have proved so
far we know that $Y_i:=f^{-1}_i(L)$ is connected and non-empty for
every $i=1,...,m$. Clearly, $Y=Y_1\cup ...\cup Y_m$. To prove that
$Y$ is connected it will be then sufficient to show that we can
reorder the components $X_1,...,X_m$ (possibly with repetitions, by
increasing $m$ if necessary) in such a way that $Y_{i-1}\cap Y_i\neq
\emptyset$ for every $i=2,...,m$. To do that we use proposition (1.2)
to get a reordering (possibly with repetitions) $X_1,...,X_m$ such
that ${\dim}(X_{i-1}\cap X_i)\geq d$ for every $i=2,...,m$. Therefore
for every $i=2,...,m$ there is an irreducible component $Z_i$ of
$X_{i-1}\cap X_i$ such that ${\dim}(Z_i)\geq d$. Since $d\geq r$,
then apply the last part of theorem (1.3) to the
$u_i:=f/Z_i:Z_i\rightarrow P$ to deduce that $u^{-1}_i(L)\neq\emptyset$
for every $i=2,...,m$. Since $u^{-1}_i(L)\subset Y_{i-1}\cap Y_i$,
it follows that $Y_{i-1}\cap Y_i\neq\emptyset$ for every $i=2,...,m$,
as desired.

Therefore proposition (2.1) is proved in case $d=r$. The case
$d>r$ can be reduced to the case $d=r$ as follows. Let $W$ be a
closed subscheme of $Y$ of dimension $<d-r$. Then $f(W)$ is a
closed subscheme of $Y$ of dimension $<d-r$. Pick a sufficiently
large integer $a>0$ such that ${\Cal O}_P(a)$ is a very ample
line bundle on $P$ such that $S_a\cong H^0(P,{\Cal O}_P(a))$.
The existence of such an $a$ comes from the fact that $S$ is a
finitely generated graded $k$-algebra. Let
$t_{r+1},...,t_d\in S_a=H^0(P,{\Cal O}_P(a))$ be $d-r$ general
homogeneous elements of degree $a$ of $S$. Since ${\dim}(f(W))<d-r$ and
$t_{r+1},...,t_d$ are general, we infer that
$f(W)\cap V_{+}(t_{r+1},...,t_d)=\emptyset$, whence $W\cap Y'=\emptyset$,
where $Y':=f^{-1}(V_{+}(t_1,...,t_d))$. By what we have already
seen before, $Y'$ is connected and non-empty. Moreover, let $Z$
be an arbitrary irreducible component of $Y$. Since $X$ is
$d$-connected and $Y=f^{-1}(V_{+}(t_1,...,t_r))$, ${\dim}(Z)>d-r$,
and hence ${\dim}(f(Z))>d-r$ because $f$ is finite. In particular,
$f(Z)$ meets $V_{+}(t_{r+1},...,t_d)$, or else, $Z$ meets $Y'$.
Suppose that $Y\setminus W$ is disconnected; then $Y'$ is also
disconnected because $Y'$ meets every irreducible component of
$Y$ and $Y'\cap W=\emptyset$, a contradiction.  Q.E.D.
\enddemo

(2.2) We shall show that a construction used by Deligne (see [3], or
also [6]) to simplify the proof of Fulton-Hansen connectedness theorem
can easily be generalized to weighted projective spaces.
Having the system $e=(e_0,e_1,...,e_n)$ of weights fixed, consider
the weighted projective space
$$ {\Bbb P}^{2n+1}(e,e)={\Proj}(k[T_0,...,T_n;U_0,...,U_n]), $$
where $T_0,...,T_n,U_0,...,U_n$ are $2n+2$ indeterninates over $k$ such
that ${\deg}(T_i)={\deg}(U_i)=e_i$ for every $i=0,1,...,n$.
Consider the closed subschemes
$$ L_1=V_{+}(T_0,...,T_n) \; \text{and}\; L_2=V_{+}(U_0,...,U_n) $$
of $P:={\Bbb P}^{2n+1}(e,e)$. Then $L_1\cap L_2=\emptyset$. Set
$U:=P\setminus (L_1\cup L_2)$. Since $T_i-U_i$ is a homogeneous
element of degree $e_i$, it makes sense to consider also the closed
subscheme $H:=V_{+}(T_0-U_0,...,T_n-U_n)$ of $P$. Clearly, $H\subset U$.
The two natural inclusions $k[T_0,...,T_n]\subset k[T_0,...,T_n;U_0,...,U_n]$
and $k[U_0,...,U_n]\subset k[T_0,...,T_n;U_0,...,U_n]$ yield two
rational maps $g_i:{\Bbb P}^{2n+1}(e,e)\dashrightarrow {\Bbb P}^{n}(e)$,
$i=1,2$, which give rise to the rational map
$$
g:{\Bbb P}^{2n+1}(e,e)\dashrightarrow {\Bbb P}^{n}(e)\times{\Bbb P}^{n}(e).
$$
Then $g$ is defined precisely in the open subset $U$ of ${\Bbb P}^{2n+1}(e,e)$.
Alternatively, if we interpret ${\Bbb P}^{n}(e)$ as the geometric
quotient $(k^{n+1}\setminus\lbrace{0}\rbrace)/{\Bbb G}_m$ (where the action
of the multiplicative group ${\Bbb G}_m=k\setminus \lbrace{0}\rbrace$
on $k^{n+1}\setminus\lbrace{0}\rbrace$ is given by
$\lambda (t_0,...,t_n):=({\lambda}^{e_0}t_0,...,{\lambda}^{e_n}t_n)$,
with $\lambda \in{\Bbb G}_m$ and $(t_0,...,t_n)\in k^{n+1}\setminus
\lbrace{0}\rbrace$), then the map $g$ is defined by
$$ g([t_0,...,t_n;u_0,...,u_n])=([t_0,...,t_n],[u_0,...,u_n]).$$
It is clear that $g/H$ defines an isomorphism $H\cong \Delta$.
Consider the commutative diagram
$$\CD
X'     @>g'>>    X \\
@Vf'VV           @VVfV \\
U     @>g>>      {\Bbb P}^{n}(e)\times{\Bbb P}^{n}(e) \\
@AAA             @AAA  \\
H      @>>>      \Delta
\endCD $$
where the top square is cartesian, the vertical arrows of the bottom
square are the canonical closed immersions, and the bottom horizontal
arrow is an isomorphism. We shall need the following fact:

(2.3) In the hypotheses of theorem (0.1), the variety
$X'=X\times_{{\Bbb P}^{n}(e)\times{\Bbb P}^{n}(e)}U$
is $(d+1)$-connected.

Indeed, since $X$ is $d$-connected, by proposition (1.2) we can
reorder the irreducible components $X_1,...,X_m$ of $X$ (possibly
with repetitions, by increasing $m$ if necessary) so that
${\dim}(X_{i-1}\cap X_i)\geq d$ for every $i=2,...,m$. Set
$X'_i:=g^{'-1}(X_i), \; i=1,...,m$. Observe that every closed fibre of
$g$ (and hence also of $g'$) is irreducible. Indeed, every closed fibre
of $g$ is in fact isomorphic to ${\Bbb G}_m$. It follows that
$X'_i$ is irreducible for every $i=1,...,m$, and in particular,
$X'_1,...,X'_m$ are the irreducible components of $X'$. Moreover,
since ${\dim}(X_{i-1}\cap X_i)\geq d$ we get that
${\dim}(X'_{i-1}\cap X'_i)\geq d+1$ for every $i=1,...,m$. This
proves the claim of (2.3) via proposition (1.2).

(2.4) Now we can easily prove the first part of theorem (0.1).
By (2.3) and proposition (2.1) (applied to the finite morphism
$f':X'\rightarrow U\subset {\Bbb P}^{2n+1}(e,e)$ and $L:=H$,
with $r=n+1\leq d+1$ we infer that $f^{'-1}(H)$ is $(d-n)$-connected.
On the other hand, since $f^{-1}(\Delta)\cong f^{'-1}(H)$ we
get the first part of theorem (0.1).

(2.5) Now we proceed to prove the second part of theorem (0.1).
I claim that one can reduce oneself to the case when $X$ is irreducible
of dimension $d+1$ such that $d\geq n$. Indeed, if $X_1,...,X_m$ are
the irreducible components of $X$ then $X_1\setminus W,...,X_m\setminus W$
are the irreducible components of $X\setminus W$. Then the claim follows
easily from corollary (1.8).

So, from now on assume $X$ irreducible of dimension $d+1$. Set $X':=f(X)$.
I claim that it is sufficient to prove the second part of theorem (0.1)
for $X'$, i.e. we may assume that $f$ is a closed embedding.

Indeed, let $W$ be a closed subscheme of $Y:=f^{-1}(\Delta)$,
of dimension $<d-r$ and set $W':=f(W)$. Then clearly
${\dim}(W)={\dim}(W')={\dim}(f^{-1}(W'))$ (because f is finite),
and $W\subset f^{-1}(W')$. Assuming that $X'\cap\Delta\setminus W'$ is
${\bold G}_3$ in $X'\setminus W'$, by theorem (1.4) we infer that
$Y\setminus f^{-1}(W')$ is ${\bold G}_3$ in $X\setminus f^{-1}(W')$.

On the other hand, since $W\subset f^{-1}(W')$ we have a commutative
diagram
$$\CD
K(X\setminus W)          @>>>  K(\widehat{X\setminus W}_{/Y\setminus W}) \\
@VVV                          @VVV  \\
K(X\setminus f^{-1}(W')) @>>>  K(\widehat{X\setminus f^{-1}(W')}_{/Y\setminus
f^{-1}(W')}) \endCD $$
in which the vertical arrows are the restriction maps. The first vertical
map is clearly an isomorphism, while the bottom horizontal map is an
isomorphism because $Y\setminus f^{-1}(W')$ is ${\bold G}_3$ in
$X\setminus f^{-1}(W')$. If we prove claim (2.6) below it will follow
that the second vertical map is injective, and hence bijective because
we just saw that the composition of it with the top horizontal map of
the above commutative diagram is an isomorphism.

(2.6) The ring $K(\widehat{X\setminus W}_{/Y\setminus W})$ is a field.

To prove this we use proposition (1.6). So, if $u:Z\rightarrow X$
is the normalization morphism it is sufficient to check that
$u^{-1}(Y\setminus W)$ is connected. This follows from the first
(i.e. the connectivity) part
of theorem (0.1) applied to the finite morphism
$f\circ u:Z\rightarrow {\Bbb P}^{n}(e)\times{\Bbb P}^{n}(e)$
and the subscheme $u^{-1}(W)$ of $Z$, because
$$ u^{-1}(Y\setminus W)=(f\circ u)^{-1}(\Delta)\setminus u^{-1}(W). $$

(2.7) Summing up, to prove the second part of theorem (0.1),
we may assume that $f$ is a closed embedding and that $X$ is
irreducible. In other words, from now on $X$ is an irreducible
closed subscheme of ${\Bbb P}^{n}(e)\times{\Bbb P}^{n}(e)$.

According to the construction and notations of (2.2), set
$U_X:=g^{-1}(X)\subset U$, and denote by $Z$ the closure of $U_X$
in ${\Bbb P}^{n}(e,e)$. Denote by $g':U_X\rightarrow X$ the
restriction of $g$ to $U_X$. Then $g'$ can be also considered as
a rational map $g':Z\dashrightarrow X$ which is defined precisely
in the open subset $U_X$ of $Z$.

Let $h:X_1\rightarrow Z$ be a proper morphism with the following
properties:

- $h$ is birational and the restriction $h/h^{-1}(U_X)$ defines an
isomorphism between $h^{-1}(U_X)$ and $U_X$;

- the composition $f:=g'\circ h:X_1\rightarrow X$ is a proper morphism
(in particular, $f$ is everywhere defined).

The existence of $(X_1,h)$ is obvious because one can take as $X_1$
the closure of the graph ${\Gamma}_{g'}\subset U_X\times X$ in
$Z\times X$, and as $h$ the restriction to $X_1$ of the projection
$Z\times X\rightarrow Z$. Then the restriction to $X_1$ of the
projection $Z\times X\rightarrow X$ is a proper surjective morphism
$f(=g'\circ h)$, so that it makes sense to speak about the field extension
$f^{\ast}:K(X)=K(X\setminus W)\rightarrow K(X_1\setminus f^{-1}(W))$.

Applying theorem (1.4) to the morphism
$f:X_1\setminus f^{-1}(W)\rightarrow X\setminus W$, we get that
$K(\widehat{X_1\setminus f^{-1}(W)}_{/f^{-1}(X\cap\Delta\setminus W})
\cong [K(\widehat{X\setminus W}_{X\cap\Delta\setminus W})\otimes_
{K(X\setminus W)} K(X_1\setminus f^{-1}(W))]_0$. Then the key point of the
proof of the second part of theorem (0.1) is the following claim:

(2.8) The canonical map
$K(X_1\setminus f^{-1}(W))\rightarrow K(\widehat{X_1\setminus f^{-1}(W)}_
{/f^{-1}(X\cap\Delta\setminus W)})$ is an isomorphism.

Accepting the claim (2.8) for the moment, we see that via the above
isomorphism, the claim obviously implies the following fact:

(2.9) The canonical map
$K(X\setminus W)\rightarrow K(\widehat{X\setminus W}_{X\cap\Delta\setminus W})$
is an isomorphism, i.e. $X\cap\Delta\setminus W$ is ${\bold G}_3$ in
$X\setminus W$. Therefore, in view of the reduction made at (2.5) we
proved the last part of theorem (0.1) modulo the claim (2.8)

(2.10) Now we proceed to the proof the claim of (2.8). Since by (2.2)
$g/H$ defines an isomorphism $H\cong\Delta$ then $g'/U_X\cap H$ defines
an isomorphisms $Z\cap H = U_X\cap H\cong X\cap\Delta$. By the
construction of $h$ the subvariety $Y:=h^{-1}(Z\cap H)$ is isomorphic
to $Z\cap H\cong X\cap\Delta$ via $f/Y$. Since
$Y\setminus f^{-1}(W)\subset f^{-1}(X\cap\Delta\setminus W)$ we get
the canonical morphisms
$$ \widehat{X_1\setminus f^{-1}(W)}_{/Y\setminus f^{-1}(W)}\rightarrow
\widehat{X_1\setminus f^{-1}(W)}_{/f^{-1}(X\cap\Delta\setminus W)}
\rightarrow X_1\setminus f^{-1}(W),$$
which yield the homomorphisms of $k$-algebras
$$ \CD
K(X_1\setminus f^{-1}(W))
@>a>> K(\widehat{X_1\setminus f^{-1}(W)}_{/f^{-1}(X\cap\Delta\setminus W)}) \\
   @. @VVbV \\
   @. K(\widehat{X_1\setminus f^{-1}(W)}_{/Y\setminus f^{-1}(W)})
\endCD $$

Now, claim (2.8) follows from the following two claims:

(2.8.1) The composition $b\circ a$ is an isomorphism.

(2.8.2) $K(\widehat{X_1\setminus f^{-1}(W)}_
{/f^{-1}(X\cap\Delta\setminus W)})$ is
a field.

Indeed, by (2.8.2) the map $b$ is injective, and whence by (2.8.1), an
isomorphism. By (2.8.1) again the map $a$ is an isomorphism.

It remains therefore to prove (2.8.1) and (2.8.2).

To prove (2.8.1) observe that since $H\subset U$ and
$h/h^{-1}(U_X):h^{-1}(U_X)\rightarrow U_X=Z\cap U$ is an isomorphism,
we get
$$K(X_1\setminus f^{-1}(W))\cong K(Z\setminus h(f^{-1}(W)\cap Y)),
\;\text{and}$$
$$K(\widehat{X_1\setminus f^{-1}(W)}_{/Y\setminus f^{-1}(W)})\cong
K(\widehat{Z\setminus h(f^{-1}(W))}_{/Z\cap H\setminus h(f^{-1}(W)\cap Y)}).
$$

So, we are reduced to prove that the canonical map
$$K(Z\setminus W_1)\rightarrow K(\widehat{Z\setminus V}_
{/Z\cap H\setminus V}) $$
is an isomorphism, where $V:=h(f^{-1}(W)\cap Y)$ is a closed subscheme
of $Z\cap H$ such that ${\dim}(V)={\dim}(W)<{\dim}(X)-n-1={\dim}(Z)-(n+1)-1$,
and ${\codim}_{{\Bbb P}^{2n+1}(e,e)}(H)=n+1$. (Recall that
$H=V_{+}(T_0-U_0,...,T_n-U_n)$.)
This will follow from the following more general assertion:

(2.10.1) For every closed irreducible subvariety $Z$ of ${\Bbb P}^{m}(e)=
{\Proj}(k[T_0,...,T_m])$ (with ${\deg}(T_i)=e_i, i=0,1,...,m$), for
every subscheme $H$ of the form $H=V_{+}(T_{i_1},...,T_{i_r})$, and
for every closed subscheme $V$ of $Z\cap H$ such that
${\dim}(V)<{\dim}(Z)-r-1$, then $Z\cap H\setminus V$ is ${\bold G}_3$
in $Z\setminus V$.

To do this, consider the usual projective space
${\Bbb P}^n={\Proj}(k[U_0,...,U_m])$, ${\deg}(U_i)=1, i=0,1,...,m$.
Consider the homomorphism of graded $k$-algebras
$$ \varphi:k[T_0,...,T_m]\rightarrow k[U_0,...,U_m] $$
defined by $\varphi (T_i)={U_i}^{e_i}$, $i=0,1,...,m$. Then the
morphism
$$ u:={\Proj}(\varphi):P':={\Bbb P}^n\rightarrow P:={\Bbb P}^{n}(e) $$
is finite and surjective.

(2.10.2) I claim that for every closed subvariety $Z$ of $P$, then
$u^{-1}(Z)_{red}$ is $d$-connected.

Indeed, consider the morphism $f:=u\times i:P'\times Z\rightarrow P\times P$,
where $i:Z\rightarrow P$ is the embedding of $Z$ in $P$. Since $Z$ is
$d$-connected, $P'\times Z$ is $(n+d)$-connected (as is easy to see).
Then by the first part of theorem (0.1) (already proved),
$f^{-1}(\Delta)\cong u^{-1}(Z)$ is $d$-connected, as claimed.

In particular, coming back to (2.10.1), since $Z$ is irreducible,
by (2.10.2), $Z':=u^{-1}(Z)$ is $({\dim}(Z)-1)$-connected. If
we set $H':=u^{-1}(H)_{red}$ then $H'$ is the linear subspace
of ${\Bbb P}^n$ given by $U_{i_1}=...=U_{i_r}=0$. Thus, if
$A:=Z\cap H$ then $u^{-1}(A)_{red}=Z'\cap H'$. By the result
of Faltings in the form of corollary (1.9) we infer that the
natural map
$$K(Z'\setminus u^{-1}(V))\rightarrow
K(\widehat{Z'\setminus u^{-1}(V)}_{/u^{-1}(A\setminus V)}) $$
is an isomorphism. On the other hand, the latter ring is by
theorem (1.4) (which can be applied because the fact that $Z'$ is
$({\dim}(Z)-1)$-connected implies that every irreducible
component of $Z'$ dominates $Z$) canonically isomorphic to
$$ [K(\widehat{Z\setminus V}_{/A\setminus V})\otimes_{K(Z\setminus V)}
K(Z'\setminus u^{-1}(V))]_0 .$$
These two isomorphisms imply (2.10.1). In this way (2.8.1) is
proved.

(2.11) Now we shall prove (2.8.2). Let $v:X'_1\rightarrow X_1$ be
the normalization of $X_1$, and set $f':=f\circ v:X'_1\rightarrow
X\subset {\Bbb P}^{n}(e)\times {\Bbb P}^{n}(e)$. Using proposition
(1.6), all we have to check is that
$v^{-1}(f^{-1}(X\cap\Delta\setminus W))=f^{'-1}(\Delta\setminus W)$
is connected. Let $X'_1 @>v'>> T @>v">> X\subset {\Bbb P}^{n}(e)
\times {\Bbb P}^{n}(e)$ be the Stein factorization of $f'$, i.e.
$f'=v"\circ v'$, where $v'$ has connected fibres, and $v"$ is a
finite morphism. By the first part of theorem (0.1) (already proved),
$v"^{-1}(\Delta\setminus W)$ is connected because
${\dim}(v"^{-1}(W))={\dim}(W)<{\dim}(X)-n-1$, by hypotheses, and hence
$v^{'-1}(v^{"-1}(\Delta\setminus W))=f^{'-1}(\Delta\setminus W)$ is
also connected because $v"$ has connected fibres. Thus (2.8.2) is
proved.

In this way the proof of theorem (0.1) is complete.

\subhead{Remark}\endsubhead It is well known that at least in
characteristic zero, ${\Bbb P}^n(e)$ appears as the quotient of ${\Bbb P}^n$
by a finite group $G$ (which is the product of the cyclic groups of order
$e_i, i=0,...,n$ (see [3])). Then one may ask whether theorem (0.1)
holds true for every quotient $P':={\Bbb P}^n/G$ of ${\Bbb P}^n$ by a
finite group $G$ of automorphisms of ${\Bbb P}^n$. The answer is no in
general. Indeed, consider the action of the group $G$ of roots of
order $n+2$ of $1$, such that $n\geq 3$, $n+2$ is prime and different from
$\text{char}(k)$, on ${\Bbb P}^{n+2}$ given by
$g(t_0,t_1,...,t_n)=(t_0,gt_1,g^{2}t_2,...,g^{n}t_n)$. Then the Fermat's
hypersurface $F$ of ${\Bbb P}^n$ given by
$t^{n+2}_0+t^{n+2}_1+...+t^{n+2}_n=0$ is
$G$-invariant and $G$ acts freely on $F$. If we take $X:=F\times F$,
and as $f:X\rightarrow P'\times P'$ the composition
$$\CD F\times F \subset {\Bbb P}^{n}\times {\Bbb P}^{n} @>{u\times u}>>
P'\times P', \endCD $$
(with $u:{\Bbb P}^{n}\rightarrow P'$ the canonical morphism), then
$f^{-1}({\Delta}_{P'})$ has $n+2$ connected components, although $X$ is
$(2n-3)$-connected ($X$ is irreducible of dimension $2n-2$)
and ${\dim}(P')=n\leq 2n-3$, if $n\geq 3$ (see the
proof of theorem (3.5) and the example following it, for details).

\head{3. Proofs of theorems (0.2) and (0.3) and other consequences}
\endhead

First we mention the following direct consequence of theorem (0.1).

\proclaim{Corollary (3.1)} Let $Y$ and $Z$ be two irreducible subvarieties
of ${\Bbb P}^{n}(e)$ such that ${\dim}(Y)+{\dim}(Z)\geq n+1$. Then
$Y\cap Z\setminus W$ is ${\bold G}_3$ in $Y\times Z\setminus W$ for
every closed subscheme $W$ of $Y\cap Z\cong (Y\times Z)\cap\Delta$ such
that ${\dim}(W)<{\dim}(Y)+{\dim}(Z)-n-1$. In particular, if $Y$ is a
closed irreducible subvariety of ${\Bbb P}^{n}(e)$ of dimension
$>\frac{n}{2}$, then ${\Delta}_Y\setminus W$ is ${\bold G}_3$ in
$Y\times Y\setminus W$ for every closed subscheme $W$ of the diagonal
${\Delta}_Y$ of $Y\times Y$ such that ${\dim}(W)<2{\dim}(Y)-n-1$.
\endproclaim

\demo{Proof} Apply theorem (0.1) to the inclusion of
$X:=Y\times Z$ in ${\Bbb P}^{n}(e)\times {\Bbb P}^{n}(e)$. Q.E.D.
\enddemo

(3.2) The last part of corollary (3.1) is just the first part of
theorem (0.2). In particular, ${\Delta}_Y$ is ${\bold G}_3$ in
$Y\times Y$ for every irreducible subvariety $Y$ of ${\Bbb P}^{n}(e)$
of dimension $>\frac{n}{2}$. This latter statement implies the second
part of theorem (0.2). The argument to do that is essentially due to
Speiser (see [20], or also [15], pp. 200-201), at least in the case
when $Y$ is smooth. Below I shall only indicate the necessary steps
in order to make Speiser's proof of the implication "${\Delta}_Y$
${\bold G}_3$ in $Y\times Y$ implies ${\Lef}(Y\times Y,{\Delta}_Y)$"
work in this more general situation.

First, exactly as in Speiser's proof (loc.cit.), we can reduce ourselves
to the case when the vector bundle $E$ is of rank one, i.e. a line
bundle. In this case, using the irreducibility of $Y$, we may assume
that $E$ is a subsheaf of the constant sheaf $K(Y\times Y)$ (see e.g.
[14], II, prop. 6.15). Then all we have to check are the following
two claims:

(3.2.1) For every $x:=(y,y)\in{\Delta}_Y$ we have
$$ {\Cal O}_{\widehat{Y\times Y},x}\cap K(X)={\Cal O}_{Y\times Y,x},$$
where the intersection is taken in the total ring of fractions of
${\Cal O}_{\widehat{Y\times Y},x}$.

Using (3.2.1), the ${\bold G}_3$ property of $(Y\times Y,{\Delta}_Y)$,
and the arguments given in [15], pp. 200-201, we infer that every
section of $H^0(\widehat{Y\times Y},\hat{E})$ comes from a
section of $E$ defined in a suitable open neighbourhood $U$ of
${\Delta}_Y$ in $Y\times Y$. This latter section of $E$ extends to the
whole $Y\times Y$ by using the ${\bold S}_2$ property of $Y$ (which
implies the ${\bold S}_2$ property of $Y\times Y$) and the following:

(3.2.2) $\text{Codim}_{Y\times Y}(Y\times Y\setminus U)\geq 2$.

We first prove (3.2.1). In the case when $Y$ is normal a proof can
be found in [20], pp. 16-17 (see also [15], p. 209, Ex. 4.12).
However, as G. Chiriacescu pointed out, (3.2.1) comes from the
following general statement. If $A\rightarrow B$ is a flat homomorphism
of local rings, then $A=B\cap [A]_0$, where the intersection is
taken in the total ring of fractions $[B]_0$ of $B$. To prove this
statement, let $b=\frac{u}{v}\in B$, with $u,v\in A$, and $v$ a
non-zero divisor of $A$. If $v$ is invertible there is nothing to
prove. If $v$ belongs to the maximal ideal of $A$ then $u=bv\in
Bv\cap A=Av$ (the latter equality holds because the map
$A\rightarrow B$ is faithfully flat). It follows that $b\in A$, as
desired.

Now we prove (3.2.2). Set $X=Y\times Y$. Since ${\dim}(Y)>\frac{n}{2}$
then $X$ is a closed subvariety of ${\Bbb P}^{n}(e)\times{\Bbb P}^{n}(e)$
of dimension $>n$. Then the claim follows from the construction of the
proof of theorem (0.1), because if $D$ were a hypersurface of $X$
which does not intersect ${\Delta}_Y=\Delta\cap X$ then (in the notations
of the proof of theorem (0.1)) $g^{'-1}(D)$ would be a hypersurface
of $U_X$ which does not intersect $U_X\cap H$ ($U_X\cap H\subset
g^{'-1}({\Delta}_Y)$), where $g'=g/U_X:U_X\rightarrow X$ is a morphism
whose closed fibres all all isomorphic to ${\Bbb G}_m$. Taking the
closure $Z$ of $U_X$ in ${\Bbb P}^{2n+1}(e,e)$, and the closure $D'$ of
$g^{'-1}(D)$ we would get a hypersurface of $Z$ which does not
intersect $Z\cap H$. But recalling that $Z$ is a closed irreducible
subvariety of ${\Bbb P}^{2n+1}(e,e)$ of dimension $\geq n+2$ and
$H=V_{+}(T_0-U_0,...,T_n-U_n)$, this fact is impossible. Thus
(3.2.2) is proved.

This finishes the proof of theorem (0.2). Q.E.D.

(3.3) Theorem (0.3) is a direct consequence of theorems (0.2) and (1.11).

(3.4) Now consider the following situation:

- $Y$ a smooth projective subvariety of ${\Bbb P}^{n}(e)$ of
dimension $>\frac{n}{2}$.

- $G$ a finite abelian group of order $d$ acting freely on $Y$ such
that $d\geq 2$, and if $\text{char}(k)>0$, then $d$ is prime to
$\text{char}(k)$.

- Denote by $X$ the quotient $Y/G$.

Then we prove:

\proclaim{Theorem (3.5)} In the situation of (3.4) we have
$$ [K(\widehat{X\times X}):K(X\times X)]=d, $$
where $K(\widehat{X\times X})$ is the field of formal rational
functions of $X\times X$ along the diagonal ${\Delta}_X$. In particular,
${\Delta}_X$ is ${\bold G}_2$ (but not ${\bold G}_3$) in $X\times X$
(in the terminology of [16], or also [15]).
\endproclaim

\demo{Proof} Denote by $u:Y\rightarrow X$ the canonical $\acute{e}$tale
morphism. Then
$$ (u\times u)^{-1}({\Delta}_X)=(G\times G){\Delta}_Y, $$
where $G\times G$ acts on $Y\times Y$ in the canonical way.

A simple computation shows that if $g,g',h,h'\in G$ then
$$ (g,gh){\Delta}_Y\cap (g',g'h'){\Delta}_Y\neq\emptyset \;\text{iff}\;
h=h'\;\text{iff}\;(g,gh){\Delta}_Y=(g',g'h'){\Delta}_Y. $$

This shows that
$$ (u\times u)^{-1}({\Delta}_X)={\Delta}_1\cup ...\cup{\Delta}_d ,
\tag 3.5.1 $$
where ${\Delta}_i=(e,h_i){\Delta}_Y, \;i=1,...,d$, with
$\lbrace h_1=e,...,h_d \rbrace $ is a full system of representatives
for the quotient group $(G\times G)/{\Delta}_G$ and $e$  the unity of $G$
(and in particular, ${\Delta}_1={\Delta}_Y$). It follows that
${\Delta}_i\cap{\Delta}_j=\emptyset$ if $i\neq j$, and that for any
$i$ and $j$ there is an element $(e,h_{ij})\in G\times G$ such
that $(e,h_{ij}){\Delta}_i={\Delta}_j$. From (3.5.1) it follows
$$ K(\widehat{Y\times Y}_{/(u\times u)^{-1}({\Delta}_X)})=
{\prod}_{i}K(\widehat{Y\times Y}_{/{\Delta}_i})\cong
\prod K(\widehat{Y\times Y}_{/{\Delta}_Y}). \tag 3.5.2 $$
(The last product of (3.5.2) has $d$ factors.) On the other hand,
by theorem (1.4) we have
$$ K(\widehat{Y\times Y}_{/(u\times u)^{-1}({\Delta}_Y)})\cong
[K(\widehat{X\times X}_{/{\Delta}_X})\otimes_{K(X\times X)}K(Y\times Y)]_0.
\tag 3.5.3 $$

By corollary (3.1) we have $K(\widehat{Y\times Y}_{/{\Delta}_Y})
\cong K(Y\times Y)$, and therefore (3.5.2) and (3.5.3) yield:
$$ K(\widehat{X\times X}_{/{\Delta}_X})\otimes_{K(X\times X)}
K(Y\times Y)\cong \prod K(Y\times Y) \;(d\;\text{times}).\tag 3.5.4 $$
Here we have used the fact that the total ring of fractions of a
product of fields coincides to that product itself. Finally, by
proposition (1.6), $K(\widehat{X\times X}_{/{\Delta}_X})$ is a
field, and hence by (3.5.4) one concludes the proof of theorem
(3.5). Q.E.D.
\enddemo

\subhead{Example}\endsubhead If in theorem (3.5) we take as $Y$
the Fermat's surface of equation $T^{5}_0+T^{5}_1+T^{5}_2+T^{5}_3=0$
with the action of the multiplicative group $G$ of roots of 1 of
order 5 ($\text{char}(k)\neq 5$) acting by
$$ g(t_0,t_1,t_2,t_3)=(t_0,gt_1,g^{2}t_2,g^{3}t_3), \; g\in G , $$
then we get that the Godeaux surface $X=Y/G$ has the property that
$K(\widehat{X\times X})$ is a field extension of degree 5 of
$K(X\times X)$. In particular, ${\Delta}_X$ is ${\bold G}_2$ (but
not ${\bold G}_3$) in $X\times X$.

\head{Appendix}\endhead

In this appendix we give a proof of theorem (1.3) of Grothendieck.
Assume first that $X$ is irreducible of dimension $n\geq r$, and
set $P:={\Proj}(S)$ and $Y:=f^{-1}(V_{+}(t_1,...,t_r))$. Since
$P\setminus V_{+}(t_1,...,t_r)=D_{+}(t_1)\cup ...\cup D_{+}(t_r)$,
$f$ is finite and $D_{+}(t_i)$ is affine, it follows that $X\setminus Y$
is the union of the affine open subsets $f^{-1}(D_{+}(t_i))$, $i=1,...,r$.
If $Z$ is an algebraic variety, consider Hartshorne's cohomological
dimension $\text{cd}(Z)$ of $Z$ defined by
$\text{cd}(Z):=\text{max}\lbrace m\geq 0/H^m(Z,F)\neq 0 \;\text{for some
coherent sheaf}\; F\in\text{Coh}(Z)\rbrace$. It follows that
$\text{cd}(X\setminus Y)\leq r-1\leq n-1$ (resp.
$\text{cd}(X\setminus Y)\leq n-2$ if $r\leq n-1$).
Now, we have the following:

\proclaim{Proposition}(Hartshorne-Lichtenbaum) In the above situation, the
inequality $cd(X\setminus Y)\leq n-1$ (resp. $cd(X\setminus Y)\leq n-2$)
implies $Y\neq\emptyset$ (resp. $Y$ connected). \endproclaim

\demo{Proof} The variety $X$ is projective because $f$ is a finite
morphism and $P$ is projective. The first assertion follows from the
so-called Hartshorne-Lichtenbaum theorem (see [11], [17], or [15]),
asserting in this case that
$\text{cd}(X\setminus Y)=n$ if and only if $Y=\emptyset$. Assume therefore
($X$ irreducible and) $\text{cd}(X\setminus Y)\leq n-2$, and $Y$
disconnected, with $Y=Y_1\cup Y_2$, $Y_1, Y_2$
non-empty closed subsets of $Y$ such that $Y_1\cap Y_2=\emptyset$.
Because $X$ is projective of dimension $n$, there is an invertible
sheaf $L$ on $X$ such that $H^n(X,L)\neq 0$ (just take any sufficiently
negative power of an ample line bundle on $X$).

In the exact sequence
$$\CD H^{n}_{Y_i}(X,L) @>>> H^n(X,L) @>>> H^n(X\setminus Y_i,L) \endCD $$
the last space is zero by Hartshorne-Lichtenbaum theorem because
$Y_i\neq\emptyset$ for $i=1,2$. It follows that $h^{n}_{Y_i}(X,L)\geq
h^n(X,L)$ for $i=1,2$. Moreover,
$H^{n}_{Y}(X,L)\cong H^{n}_{Y_1}(X,L)\oplus H^{n}_{Y_2}(X,L)$, and
hence, $h^{n}_{Y}(X,L)\geq 2h^n(X,L)$.

On the other hand, in the exact sequence
$$\CD
H^{n-1}(X\setminus Y,L) @>>> H^{n}_{Y}(X,L) @>>> H^n(X,L) @>>>
H^n(X\setminus Y,L)  \endCD $$
the extreme spaces are zero by $\text{cd}(X\setminus Y)\leq d-2$, whence
$h^{n}_{Y}(X,L)=h^n(X,L)$, contradicting the previous inequality
because $h^n(X,L)\geq 1$. Q.E.D.
\enddemo

\remark{Remark} Consider the Serre map of $P={\Proj}(S)$
$$\alpha :S\longrightarrow \bigoplus_{a\geq 0}H^0(P,{\Cal O}_P(a)), $$
and for $t_i\in S_{a_i}, i=1,...,r$, set
$s_i:=\alpha (t_i)\in H^0(P,{\Cal O}_P(a_i))$.
Then for a point $p\in P$ the following statement holds
(even if $S$ is not generated by $S_1$):
$p\in D_{+}(t_i)$ if and only if $s_i(p)\neq 0$.(Recall that $s_i(p)\neq 0$
means that $(s_i)_p\notin m_p{\Cal O}_P(a_i)_p$, where $m_p$ is the
maximal ideal of ${\Cal O}_p$.)  With this observation,
the proof of the above lemma yields in fact the following more general
assertion. Assume that for every $i=1,...,r$ we are given a finitely
generated graded $k$-algebra $S^i$, a homogeneous element
$t_i\in (S^i)_{+}$, and a finite morphism $f_i:X\rightarrow {\Proj}(S^i)$.
Set $s_i:=f^{\ast}_i({\alpha}_i(t_i))$ (where ${\alpha}_i$ is the Serre
map of $P_i:={\Proj}(S^i)$). Then the zero locus of
the section $(s_1,...,s_r)$ of
$\bigoplus_{i}H^0(X,f_i^{\ast}({\Cal O}_{P_i}({\deg}(t_i)))$
is $(d-r)$-connected if $X$ is $d$-connected and $d\geq r$.
\endremark
Let us proceed now to the proof of theorem (1.3). Let $X_1,...,X_m$ be
the irreducible components of $X$. Since $X$ is $d$-connected,
${\dim}(X_i)\geq d+1$ for every $i=1,...,m$. But by what we have proved
before $Y_i:=Y\cap X_i$ is connected (and non-empty) for every
$i=1,...,m$, where, as above, $Y:=f^{-1}(V_{+}(t_1,...,t_r))$. Since
$X$ is $d$-connected, by proposition (1.2), we can reorder
the components $X_1,...,X_m$ (possibly with repetitions) so that
for every $i=2,...,m$ there is an irreducible component $Z_i$ of
$X_{i-1}\cap X_i$ of dimension $\geq d$. Applying the first
part of the above lemma to $Z_i$ we infer that $Y_{i-1}\cap Z_i=
Y_i\cap Z_i=Y\cap Z_i$ is not empty. This implies that
$Y_{i-1}\cap Y_i\neq \emptyset$ for every $i=2,...,m$. In particular,
$Y$ is connected.

We prove now that $Y$ is $(d-r)$-connected if $X$ is $d$-connected.
Let $X\subset{\Bbb P}^N$ be an arbitrary projective embedding of
$X$, and let $A$ be a general linear subspace of ${\Bbb P}^N$ of
dimension $N+r-d$. Because every irreducible component of $X$ is
of dimension $\geq d+1$ and because $Y$ is locally given by $r$
equations in $X$, every irreducible component $Z$ of $Y$ is of
dimension $\geq d+1-r$. It follows that
${\dim}(Z\cap A)={\dim}(Z)+{\dim}(A)-N\geq (d+1-r)+(N+r-d)-N=1$, and
in particular, $A$ meets every irreducible component of $Y$.

Set $Y':=Y\cap A$. If $A$ is defined by linear equations
$t_{r+1}=...t_d=0$ in ${\Bbb P}^N$ then $Y'$ is just the zero locus of
the section $(t_{r+1}/X,...,t_d/X,s_1,...,s_r)$ of the sheaf
$(d-r){\Cal O}_X(1)\oplus f^{\ast}({\Cal O}_P(a_1))\oplus ...\oplus
f^{\ast}({\Cal O}_P(a_r))$, where for $i=1,...,r$,  $t_i\in S_{a_i}$,
$s_i:=f^{\ast}(\alpha (t_i))$,and
$\alpha : S\rightarrow \bigoplus_{j}H^0(P,{\Cal O}_P(j))$ is the
Serre map of $P={\Proj}(S)$. By the above lemma and the remark following
it, $Y'$ is connected.

Now, assume that $Y$ is not $(d-r)$-connected, i.e. there is a
closed subscheme $W$ of $Y$ of dimension $<d-r$ such that $Y\setminus W$
is disconnected. Since $A$ is general,
${\dim}(W\cap A)={\dim}(W)+{\dim}(A)-N<(d-r)+(N+r-d)-N=0$, or else,
$A$ does not meet $W$. Moreover, since $A$ meets every irreducible
component of $Y$, from the fact that $Y\setminus W$ is disconnected,
it follows that $Y'=Y\cap A=(Y\setminus W)\cap A$ is also disconnected,
a contradiction.  Q.E.D.

\enddocument